\documentclass[11pt,twoside, final]{amsart}
\copyrightinfo{0}{Iranian Mathematical Society}
\pagespan{1}{\pageref*{LastPage}}
\usepackage{etoolbox,lastpage}
\commby{}
\date{\scriptsize   Received: , Accepted: .}
\usepackage{amsmath,amsthm,amscd,amsfonts,amssymb,enumerate}
\usepackage{graphicx}		
\usepackage{color}
\usepackage[colorlinks]{hyperref}

\newcommand{\supp}{\mathop{\rm supp}\nolimits}
\newcommand{\skp}[2]{\left<#1,#2\right>}

\newcommand{\s}{\mathrm{sp}}

\newtheorem{theorem}{Theorem}[section]

\newtheorem{lemma}[theorem]{Lemma}
\newtheorem{corollary}[theorem]{Corollary}
\theoremstyle{definition}

\theoremstyle{remark}

\numberwithin{equation}{section}

\begin{document}


\title[Short title]{Banach-Alaoglu theorem for  Hilbert $H^*$-module}

\author[F. Author]{Zlatko Lazovi\'c}
\address[Zlatko Lazovi\'c]{Faculty of Mathematics\\University of Belgrade\\ 11 000 Belgrade\\ Serbia}
\email{zlatkol@matf.bg.ac.rs}

\thanks{$^*$Corresponding author}
%

 \maketitle
%

\begin{abstract}
We provided an analogue Banach-Alaoglu theorem for  Hilbert $H^*$-module. We construct a $\Lambda$-weak$^*$ topology on a  Hilbert $H^*$-module
over a proper  $H^*$-algebra $\Lambda$, such that the unit ball  is compact with respect to $\Lambda$-weak$^*$ topology.
\end{abstract}

\section{\bf Introduction}

\noindent A Hilbert $H^*$-module $L$ over an $H^*$-algebra $\Lambda$ is a right $\Lambda$-module which possesses
a $\tau(\Lambda)$-valued product, where  $\tau(\Lambda)=\{ab|\,\, a,b\in \Lambda\}$ is the trace-class. At the same time, $H$ is a Hilbert space with the inner
product given by the action of the trace on the $\tau(\Lambda)$-valued product.

The notion of $H^*$-module is introduced by Saworotnow in \cite{Saworotnow1} under the name of generalized
Hilbert space. It has been studied by Smith \cite{Smith}, Giellis \cite{Giellis}
Molnar \cite{Molnar}, Cabrera \cite{Cabrera}, Martinez
 et al. \cite{Cabrera}, Baki\'{c} and Gulja\v{s} \cite{Damir} and others. Saworotnow has proved that the trace-class in a proper $H^*$-algebra $\Lambda$ has pre-dual. For Hilbert $H^*$-modules,  a generalization of Reisz theorem holds, i.e. for each bounded $\Lambda$-linear functional on $L$, there is $x_f\in L$ such that $f(x)=[x_f,x]$ for all $x\in L$.

Paschke \cite{Paschke} showed that self-dual Hilbert $W^*$-modules are dual Banach spaces and found topology on module such that the unit ball is compact.

In the present paper we find a topology on Hilbert $H^*$-module $L$ over a proper $H^*$-algebra $\Lambda$ such that the unit ball in $L$ is compact with respect to this topology.

\section{\bf Basic notations and preliminary results}
\label{sec:1}

\noindent We recall that an $H^*$-algebra is a complex associative Banach
algebra with an inner product $\skp{\cdot}{\cdot}$ such that
$\skp{a}{a}=||a||^2$ for all $a\in \Lambda$, and for each $a\in
\Lambda$ there exists some $a^*\in \Lambda$ such that
$\skp{ab}{c}=\skp{b}{a^*c}$ and $\skp{ba}{c}=\skp{b}{ca^*}$ for
all $b,c\in \Lambda$. The adjoint $a^*$ of $a$ need not be unique
(see \cite{Warren Ambrose}). Proper $H^*$-algebra $\Lambda$ is an
$H^*$-algebra which satisfies $a\Lambda=0 \Rightarrow a=0$ (or
$\Lambda a=0 \Rightarrow a=0$). An $H^*$-algebra $\Lambda$ is proper  if and only
if each  $a\in \Lambda$ has a unique adjoint $a^*\in \Lambda$ (see
\cite[Theorem 2.1]{Warren Ambrose}). An $H^*$-algebra is simple algebra if it has no nontrivial closed
two-sided ideals.

The trace-class in a proper $H^*$-algebra   $\Lambda$ is defined
as the set $\tau(\Lambda)=\{ab|\,\, a,b\in \Lambda\}$. The
trace-class is selfadjoint ideal of $\Lambda$ and it is dense in
$\Lambda$, with norm $\tau(\cdot)$. The norm $\tau$ is related to
the given  norm $||\cdot||$ on $\Lambda$ by $\tau(a^*a)=||a||^2$
for all $a\in\Lambda$. There exists a continuous linear form
 $\mathrm{sp}$ on $\tau(\Lambda)$
(trace) satisfying $\mathrm{sp}(ab)=\mathrm{sp}(ba)=\skp{a^*}{b}.$
In particular,
$\mathrm{sp}(a^*a)=\mathrm{sp}(aa^*)=\skp{a}{a}=||a||^2=\tau(a^*a)$.

 A Hilbert $\Lambda$-module is a right module $L$ over a
proper $H^*$-algebra $\Lambda$ provided with a mapping
$[\cdot,\cdot]\colon L\times L\rightarrow \tau(\Lambda)$, which satisfies the
following conditions: $[x,\alpha y]=\alpha[x,y];$ $[x,y+z]=[x,y]+[x,z];$ $[x,ya]=[x,y]\,a;$ $[x,y]^*=[y,x]$; $L$ is a Hilbert space with the inner product
$\skp{x}{y}=\mathrm{sp}\left([x,y]\right)$ for all $\alpha\in
\mathbb{C},$  $x,y,z\in L,$ $a\in\Lambda$ and for all $ x\in L,$ $x\neq0$ there is $a\in\Lambda,$ $a\neq0$ such that $[x,x]=a^*a.$

A $\Lambda$-linear functional on $L$ is a mapping $w'\colon L\rightarrow\tau(\Lambda)$ such that $w'(xa+yb)=w'(x)a+w'(y)b$ for all $x,y\in L,$ $a,b\in\Lambda.$ It is bounded if there exists $c>0$ such that $\tau(w'(x))\leq c||x||$ for all $x\in L$. In this case we define $||w'||=\supp\{c|\,\, \tau(w'(x))\leq c||x||\,\,\text{for all}\,\,x\in L\}$. The norm space of all bounded $\Lambda$-linear functional $f$ on $L$ we denoted by $L'$.

Let $R(\Lambda)$ be the set of all bounded linear operators
$S$ on $\Lambda$ such that $S(xy)=(Sx)y$ for all $x,y\in\Lambda$ and let $C(\Lambda)$ be the
closed subspace of $R(\Lambda)$ generated by the operators of the form
$La\colon x\rightarrow ax,$ $a\in\Lambda.$

We now state some theorems which will be necessary for the proof of main results.

\begin{theorem}\label{samodualnost H modula}\cite[Theorem 3]{Saworotnow1}
Each bounded $\Lambda$-linear functional $f$ on $L$ is of the form $[x_f,\cdot]$ for
some $x_f\in L$.
\end{theorem}

\begin{theorem}\label{funkcional na C(Lambda)}\cite[Lemma 1]{Saworotnow2} If  $a\in\Lambda$ then the mapping $f_a\colon S\rightarrow\mathrm{sp}(Sa)$, defined on $C(\Lambda)$, is a bounded linear functional and $||f_a||=\tau(a)$.
\end{theorem}

\begin{theorem}\label{predual H algebre}\cite[Theorem
1]{Saworotnow2} Each bounded linear functional on $C(\Lambda)$ is of the form
$f_a$ for some $a\in\tau(\Lambda)$. The correspondence $a\leftrightarrow f_a$ is an isometric isomorphism
between $\tau(\Lambda)$ and $C(\Lambda)^*$.  Also,  $\tau(\Lambda)$ is a Banach algebra.
\end{theorem}

For more details, we refer to  \cite{Saworotnow1, Saworotnow2, Warren Ambrose, Damir, Dijana}.

\section{\bf Results}
\noindent Let  $L$ be a Hilbert $H^*$-module over a proper $H^*$-algebra $\Lambda$ and let $B_1(L)$ is the unit ball in $L.$ We construct a topology on $W'$
 such that the unit ball in $L'$ is compact with respect to this topology. Define $\Lambda'$-weak$^*$ topology  on $L'$ with the base
\[L_{w'_0,x_1,..,x_n,S_1,..,S_n,\delta}=\left\{w'\in L'|\,\left|\mathrm{sp}(S_j(w'(x_j)-w'_0(x_j)))\right|<\delta,\, j=1,..,n\right\},\]
for $w'_0\in L',$ $x_j\in L,$ $S_j\in C(\Lambda),$ and $\delta>0.$

The main result of this paper, the compactness of the unit ball,  will be proven in Theorem \ref{kompaktnost jedinichne lopte u L'} and its corollary. Before that, we state and prove a useful lemma.

\begin{lemma}\label{normalnost funkcionala}
Let $\Lambda$ be a $H^*$-algebra and let $S\in C(\Lambda)$. Then the operator $S_a\colon \Lambda\rightarrow \Lambda,$ $S_a(x)=S(ax)$ belongs to $C(\Lambda)$.
\begin{proof} From
\[S_a(\lambda x)=S(a\lambda x)=\lambda S(ax)=\lambda S_a(x),\,\,S_a(xy)=S(axy)=S(ax)y=S_a(x)y\]
and
\[||S_a(x)||=||S(ax)||\leq ||S||\cdot||ax||\leq||S||\cdot||a||\cdot||x||,\]
it follows that  $S_a$ belongs to $R(\Lambda)$. If $L_{a_n}$ converges to $S$, then from
\begin{align*}
||(L_{a_na}-S_a)(x)||&=||(L_{a_na}(x)-S_a(x)||=||a_nax-S(ax)||\\
&=||L_{a_n}(ax)-S(ax)||\leq||L_{a_n}-S||\cdot||ax||\\
&\leq||L_{a_n}-S||\cdot||a||\cdot||x||,
\end{align*}
it follows that  $L_{a_na}$ converges to $S_a$. Thus $S_a\in C(\Lambda)$.
\end{proof}
\end{lemma}

\begin{theorem}\label{kompaktnost jedinichne lopte u L'} The set
\[K=\left\{w'\in L'|\,\, ||w'(x)||\leq1\,\, \text{for all}\,\, x\in B_1(L)\right\}\]
is compact in $\Lambda'$-weak$^*$ topology.
\begin{proof} The neighborhood $B_1(L)$ is absorbing because for
each $x\in L,$ $x\neq0$, there exists  a number $\gamma(x)=||x||>0$ such that $x\in \gamma(x)B_1(L)$. For all $w'\in K$, it  holds $w'(x)\leq1,$ $x\in B_1(L)$,
hence $|w'(x)|\leq||x||$, $x\in L$.

 According to Banach-Alaoglu theorem and Theorem \ref{predual H algebre}, the set $D_x=\{a\in\Lambda| ||a||\leq\gamma(x)\}$ is compact in weak$^*$-topology on $\tau(\Lambda)$ given by seminorms $p_{S}(\cdot)=|\s(S(\cdot))|,$ $S\in C(\Lambda)$. Let $\tau_1$ be the
product weak$^*$-topology  on $P=\prod_{x}D_{x}$, the cartesian product of all $D_x$, one for each
$x\in L$. Since each $D_x$ is weak$^*$-compact, it follows, from Tychonoff's theorem, that  $P$ is $\tau_1$-compact.

From definition of $P$ we have
\[P=\{f\colon L\rightarrow\tau(\Lambda)|\,\, ||f(x)||\leq||x|| \,\,\text{for all}\,\, x\in L\}.\] The set $P$ can contain $\Lambda$-nonlinear functionals.

It is clear that $K\subset L'\cap P$. It follows that $K$ inherits two topologies:
one from $L'$ (its $\Lambda'$-weak$^*$ topology, to which the conclusion of the
theorem refers) and the other, $\tau_1$, from $P$. We will see that these two topologies coincide on $K$, and that $K$ is a closed subset of $P$.

We now prove that  topologies $\tau_1$ and $\Lambda'$-weak$^*$  coincide on $K$. Fix some  $w'_0\in K$. Then
\begin{align*}
W_1=\left\{\right.w'\in L'|\,\,& \left|\s(S_j(w'(x_j)))-\s(S_j(w'_0(x_j)))\right|<\delta,\\
&\hspace{5mm}x_1,x_2,...,x_n\in L,\,S_1,...,S_n\in \left.C(\Lambda)\right\}
\end{align*}
and
\begin{align*}
W_2=\left\{f\in P\right.|\,\,& \left|\s(S_j(f(x_j)))-\s(S_j(\Lambda_0(x_j)))\right|<\delta,\\
&\hspace{5mm}x_1,x_2,...,x_n\in L,\,S_1,...,\left.S_n\in C(\Lambda)\right\}
\end{align*}
are local  bases in $(L', \Lambda\text{-weak}^*)$ and $(P,\tau_1)$, respectively. From $K\subset L'\cap P$ we have $K\cap W_1=K\cap W_2$, so topologies coincide on $K$.

Suppose $w'_0$ is in the $\tau_1$-closure of $K$. If $S$ is from $C(\Lambda)$, then, from  Lemma \ref{normalnost funkcionala}, operator  $S(a\cdot)\colon\Lambda\rightarrow \Lambda$ belongs to $C(\Lambda)$ for all $a\in\Lambda$.  For any $a,b\in\Lambda,$ $S\in C(\Lambda),$ $x,y\in L,$ $\varepsilon>0$, there is $\Lambda$-linear $w'\in K$ from $V_{xa+yb,w'_0,S,\varepsilon}\cap V_{x,w'_0,S(a\cdot),\varepsilon}\cap V_{y,w'_0,S(b\cdot),\varepsilon}$ ($\tau_1$-neighborhood of $w'_0$). Therefore, it holds
\[\left|\s(S[a(w'(x)-w'_0(x))])\right|<\varepsilon,\quad
\left|\s(S[b(w'(y)-w'_0(y))])\right|<\varepsilon,\]
\[\left|\s(S[(w'(xa+yb)-w'_0(xa+yb))])\right|<\varepsilon.\]

We have
\begin{align*}
&|\s(S(w'_0(xa+yb)-w'_0(x)a-w'_0(y)b))|\\
&=\left|\s\right.(S[(w'_0(xa+yb)-w'(xa+yb))-(w'_0(x)a-w'(x)a)\\
&\hspace{5mm}-(w'_0(y)b-w'(y)b)]\left.)\right|\\
&=\left|\s(\right.S(w'_0(xa+yb)-w'(xa+yb))-S[(w'_0(x)-w'(x))a]\\
&\hspace{5mm}-S[(w'_0(y)-w'(y))b]\left.)\right|\\
&\leq\left|\s(S(w'_0(xa+yb)-w'(xa+yb)))\right|+\left|\s(S[(w'_0(x)-w'(x))a])\right|\\
&\hspace{5mm}+\left|\s(S[(w'_0(y)-w'(y))b])\right|\\
&=\left|\s(S(w'_0(xa+yb)-w'(xa+yb)))\right|+\left|\s(aS[(w'_0(x)-w'(x))])\right|\\
&\hspace{5mm}+\left|\s(bS[(w'_0(y)-w'(y))])\right|\\
&\leq\varepsilon+\varepsilon+\varepsilon=3\varepsilon.
\end{align*}

Since $\varepsilon>0$ was arbitrary, we see that \[\s(S(w'_0(xa+yb)))=\s(S(w'_0(x)a)+S(w'_0(y)b))\]
 for all $S\in C(\Lambda),$ i.e.
  \[\s(S(w'_0(xa+yb)-w'_0(x)a-w'_0(y)b))=0\]
 for all $S\in C(\Lambda).$ For $S=L_{w'_0(xa+yb)-w'_0(x)a-w'_0(y)b)}$ we have
 \[\s((w'_0(xa+yb)-w'_0(x)a-w'_0(y)b))^*(w'_0(xa+yb)-w'_0(x)a-w'_0(y)b))=0.\]
 Hence $w'_0(xa+yb)=w'_0(x)a+w'_0(y)b$ for all $x,y\in L,$ $a,b\in\Lambda$, so $w'_0$ is $\Lambda$-linear.

 Let $x\in B_1(L)$. For arbitrary $\varepsilon>0$ and $S\in C_{\Lambda}$  there is $w'\in B_1(L')$ such that $|\s(S(w'_0(x)))-\s(S(w'(x)))|<\varepsilon.$ Hence \[|\s(S(w'_0(x)))|<|\s(S(w'(x)))|+\varepsilon\leq||S||\cdot||w'(x)||+\varepsilon\leq||S||+\varepsilon.\] Next, from Theorem \ref{funkcional na C(Lambda)} we have
 \begin{align*}
 ||w'_0(x)||&=\tau(w'_0(x))=||f_{w'_0(x)}||=\sup\limits_{S\in C(\Lambda),\,||S||=1}|\s(S(w'_0(x)))|\\
 &\leq \sup\limits_{S\in C(\Lambda),\,||S||=1} ||S||+\varepsilon\leq1+\varepsilon
  \end{align*}
  for arbitrary $\varepsilon>0$. Thus $||w'_0(x)||\leq1$.

 We have proven that $w'_0\in B_1(L')$, and that $B_1(L')$ is a closed subset of $P$.

 Since $P$ is $\tau_1$-compact, $B_1(L')$ is a closed subset of $P$, and $\tau_1$ and $\Lambda'$-weak$^*$ topology coincide on $B_1(L')$, we have that $B_1(L')$ is $\Lambda'$-weak$^*$ compact.
\end{proof}
\end{theorem}

\begin{corollary} The unit ball $B_1(L)$ in $L$ is compact in $\Lambda$-weak$^*$ topology with the base
\[W_{x_0,y_1,...,y_n,S_1,...,S_n,\delta}=\left\{x\in L|\,\, \left|\s(S_j([x,y_j]-[x_0,y_j]))\right|<\delta, j=1,...,n\right\},\]
for $x_0,y_j\in L,$ $S_j\in C(\Lambda),$ $\delta>0.$
\begin{proof} For each $w\in L'$ there is $y\in L$ such that $w'(x)=[y,x]$ for all $x\in L$ (Theorem \ref{samodualnost H modula}), so from Theorem \ref{kompaktnost jedinichne lopte u L'} it follows that the unit ball $B_1(L)$ in $L$ is compact in given topology.
\end{proof}
\end{corollary}

Acknowledgement. The author was supported in part by the Ministry of education and science, Republic of Serbia, Grant 174034.


\begin{thebibliography}{20}

\bibitem{Warren Ambrose} W. Ambrose,
\newblock Structure theorems for a special class of Banach algebras,
\newblock {\em  Trans. Amer. Math. Soc.}
\textbf{57} 
(1945), 
364--386.


\bibitem{Damir}
 D. Baki\'c and  B. Gulja\v{s},
\newblock Operators on Hilbert $H^*$-modules,
\newblock {\em J. Operator theory}
\textbf{46} 
(2001), 
123--137.

\bibitem{Cabrera} M. Cabrera, J. Martinez and A. Rodr\'iguez,
\newblock Hilbert modules revisited: Orthonormal bases and Hilbert-Schmidt operators,
\newblock {\em Glasgow Math. J}
\textbf{37} 
(1995), 
45--54.


\bibitem{Giellis} G. R. Giellis,
\newblock A characterization of Hilbert modules,
\newblock {\em Proc. Amer. Math. Soc.}
\textbf{36} 
(1972), 
440--442.

\bibitem{Dijana} D. Ili\v{s}evi\'c,
\newblock On redundance of one of the axioms of generalized normed space,
\newblock {\em Glasnik  Matemati\v{c}ki}
\textbf{37} 
(2002), 
135--141.

\bibitem{Molnar} L. Moln\'{a}r,
\newblock Modular bases in a Hilbert $A$-module,
\newblock {\em Czechoslovak Math.~J.}
\textbf{42} 
(1992), 
649--656.

\bibitem{Paschke} W.~L.~Paschke,
\newblock Inner product modules over $B^*$-algebras,
\newblock {\em Trans.~Amer.~Math.~Soc.}
\textbf{182} 
(1973), 
443--468.

\bibitem{Saworotnow1} P. P. Saworotnow,
\newblock A generalized Hilbert space,
\newblock {\em Duke Math. J.}
\textbf{35} 
(1968), 
191--197.

\bibitem{Saworotnow2}
 P. P. Saworotnow,
\newblock Trace class and centralizers of an $H^*$-algebra,
\newblock {\em Proceedings of the American Mathematical Society}
\textbf{26}, no. 1 
(Sep.,~ 1970),  
101--104.

\bibitem{Smith} J. F. Smith,
\newblock The structure of Hilbert modules,
\newblock {\em J. London Math. Soc.}
\textbf{8} 
(1975), 
741--749.

\end{thebibliography}
\end{document}